\documentclass{ifacconf}

\usepackage{graphicx}
\usepackage{xcolor}
\usepackage{natbib}        
\usepackage{amsmath,amssymb,amsfonts,mathtools,bm}

\newcounter{remarkcounter}
\newenvironment{remark}{\noindent {\bf \break Remark \refstepcounter{remarkcounter}\arabic{remarkcounter}. }\it}{

}

\begin{document}
\begin{frontmatter}

\title{Observer-Based State Feedback Controller for a Mindlin Plate Model in port-Hamiltonian framework \thanksref{footnoteinfo}} 

\thanks[footnoteinfo]{This project has received funding from the European Union's Horizon Europe research and innovative programme under the Marie Skłodowska-Curie Actions (MSCA) grant agreement No. 101073558 (ModConFlex).}

\author[First]{Ignacio Diaz Alastuey} 
\author[First]{Yann Le Gorrec} 
\author[First]{Yongxin Wu}

\address[First]{Université Marie et Louis Pasteur, SUPMICROTECH, CNRS, institut FEMTO-ST, F-25000 Besançon, France \\(e-mail: ignacio.diaz; yann.le.gorrec; yongxin.wu@femto-st.fr.)}

\begin{abstract}                
This paper generalises an early lumped observer‐based state‐feedback (OBSF) control design methodology, originally developed for one‐dimensional (1-D) boundary‐controlled port‐Hamiltonian systems, to a two‐dimensional (2-D) boundary‐controlled Mindlin plate. To this end, the 2-D port‐Hamiltonian Mindlin plate model is first introduced and then discretized using a structure‐preserving finite‐difference method on staggered grids. A controllability decomposition is subsequently applied to identify the controllable modes of the discretized model. Furthermore, the state-feedback and observer gains are designed so that the OBSF controller is strictly positive real. This guarantees the stability of the closed-loop system when the finite-dimensional OBSF controller is interconnected with the 2-D boundary-controlled Mindlin plate. Numerical simulations are finally presented to illustrate the effectiveness of the proposed method.
\end{abstract}


\end{frontmatter}

\section{Introduction}
Since its introduction in \citep{maschke1992}, the port-Hamiltonian (PH) formalism has attracted increasing attention within the systems and control community. By emphasizing energy interactions across multi-physical domains, the framework provides a unified perspective for modelling and control. It accommodates a broad class of systems, from lumped to distributed parameter systems \citep[Chapter 6]{book_vanderSchaft2017_PHS}, \citep{rashad2020}, while capturing the passivity and modularity in their structure. These properties have led to notable contributions in both theoretical analysis \citep{legorrec2004, legorrec2005} and practical control applications \citep{rodriguez2001}.

One class of distributed parameter systems addressed by PH systems is the boundary control systems (BCS), in which partial differential equations governing the system dynamics are actuated and observed only at the boundaries. For example, this applies to string models, beam theory, and transmission line models \citep{jacob2012}. Within this framework, observer-based control constitutes an effective strategy for controller design. In this context, \citep{toledo2020} proposes an early-lumping approach for the design of observer-based state feedback (OBSF) controllers that guarantees the passivity and asymptotic stability of the closed-loop system. A key advantage of this strategy is its ability to mitigate potential instabilities arising from neglected dynamics, such as the spillover effect \citep{bontsema1988}.

The main contribution of this paper is the application of the OBSF controller design strategy presented in \citep{toledo2020} to a two-dimensional boundary controlled PH system \citep{skrepek2021}. Specifically, we consider a Mindlin plate model and, by using the energy-preserving discretization strategy presented in \citep{Diaz_ECC2026}, show that this strategy can be adapted to this setting. This paper is organized as follows: Section 2 introduces the PH model of the Mindlin Plate together with its discretization. Section 3 recalls the strategy presented in \citep{toledo2020} and adapts it to the discretized PH Mindlin Plate. Section 4 presents some numerical implementation of the strategy adapted for the 2D model. Finally, Section 5 closes with some conclusions and perspectives on the following work.

\section{Mindlin Plate Model}
This section consists of two parts. The first part presents the Mindlin plate model in a port-Hamiltonian formulation, while the second part introduces the problem of interest and the spatial discretization of the model, leading to a finite-dimensional representation.
\subsection{Mindlin Plate in PH structure}
The Mindlin Plate model adopted in this work is based on the following displacement assumptions:
\begin{itemize}
	\item The thickness of the plate remains constant during deformation.
	\item The normal stress through the thickness is neglected.
	\item The mid-surface (neutral plane) of the plate undergoes a transverse displacement denoted by $u$.
	\item The orthogonal lines to the mid-surface remain straight lines after deformation, and the angles between these lines and the normal to the mid-surface are denoted by $\phi_i$.
\end{itemize}
\begin{figure}[ht]
	\centering
	\includegraphics[width=1\linewidth]{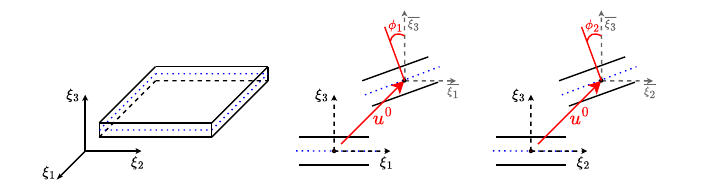}
	\caption{Mindlin Plate deformations}
	\label{fig:mindlinPlateDefor}
\end{figure}
The deformation assumptions are illustrated in Figure \ref{fig:mindlinPlateDefor}. By taking a small angle approximation given by $\tan{\phi_i}\approx\phi_i$ and $\tan^2{\phi_i}\approx0$, the displacement of an arbitrary point on the plate can be expressed as
\begin{equation}
	x=\begin{pmatrix}
		\xi_1+u_1-\xi_3\phi_1\\
		\xi_2+u_2-\xi_3\phi_2\\
		\xi_3+u_3
	\end{pmatrix}. \label{ec:MindlinPlatePosition}
\end{equation}
where $u_1,u_2,u_3,\phi_1,\phi_2\in C^1(\mathbb{R}^+_0,\Omega)$ with $\Omega$ denoting the spatial domain of the undeformed mid-surface. With these definitions, the kinetic energy of the system can be expressed in terms of the generalized momentum variables $p$ and the generalized velocity variables $M_p^{-1} p$ as
\begin{equation*}
	K(t)=\frac{1}{2}\iint_{\Omega}p^TM_p^{-1}pdA,
\end{equation*}
where 
\[
\begin{aligned}
	p&=\begin{pmatrix}
		\rho h\dot{u}_1 & \rho h\dot{u}_2 & \rho h\dot{u}_3 & \rho\tfrac{h^3}{12}\dot{\phi}_1 & \rho\tfrac{h^3}{12}\dot{\phi}_2
	\end{pmatrix}^T,\\
	M_p&=\text{diag}\left\{\begin{pmatrix}
		\rho h & \rho h & \rho h & \rho\tfrac{h^3}{12} & \rho\tfrac{h^3}{12}
	\end{pmatrix}\right\}.
\end{aligned}
\]
To formulate the elastic potential energy of the plate, the strain tensor is needed. The deformation is defined by $\mathbf{u}=x-x_0$, where $x_0$ denotes the position of an arbitrary point in the undeformed configuration and $x$ is given by \eqref{ec:MindlinPlatePosition}. Under the assumption of small deformations,  the material displacement gradient tensor is given by
\[
\boldsymbol{\epsilon}=\frac{1}{2}\Big(\nabla_\xi \mathbf{u} +(\nabla_\xi \mathbf{u})^T\Big).
\]
Using the generalized Hooke's law, which establishes a linear relation between the strain and stress tensors, the elastic potential energy of the plate can be expressed in terms of a generalized displacement coordinate $q$ as
\begin{equation*}
	U(t) = \frac{1}{2}\iint_{\Omega} q^TC_qq dA,
\end{equation*}
where 
\[\begin{array}{lccccc}
	q&=(
		\partial_1 u_1 & \partial_2 u_2 & \partial_1 u_2+\partial_2 u_1 & \partial_1 \phi_1 & \partial_2 \phi_2 \\ & \dots & \partial_1 \phi_2+\partial_2 \phi_1 & \partial_1 w^0-\phi_1 & \partial_2 w^0-\phi_2
	)^T,
\end{array}\]
\[
C_q = \frac{E}{\nu^2-1}\begin{pmatrix}
	h\mathcal{C}_q & 0_3 & 0_{3\times2}\\
	0_3 & \tfrac{h^3}{12}\mathcal{C}_q & 0_{3\times2}\\
	0_{2\times3} & 0_{2\times3} & \kappa h\frac{\nu-1}{2}I_2
\end{pmatrix},
\]
\[
\mathcal{C}_q = \begin{pmatrix}
	-1 & -\nu & 0\\ -\nu & -1 & 0\\ 0 & 0 & \frac{\nu-1}{2}
\end{pmatrix}.
\]
Here, $\nu$ denotes the Poisson's ratio, $E$ is the Young's modulus, $0_{i\times j}$ the zeros matrix of dimension $i\times j$, $0_3$ the square zero matrix of dimension $3$, $I_2$ the identity matrix of dimension $2$ and $\kappa$ is a correction factor as shear stress is known to be parabolic and not linear. Following Hamilton’s principle–based modeling approach presented in \citep{ponce2024}, the system can be expressed in port-Hamiltonian form as
\begin{equation}
	\begin{pmatrix}
		\dot{p} \\
		\dot{q}
	\end{pmatrix} = \begin{pmatrix}
			0_5 & \mathcal{J} \\ -\mathcal{J}^* & 0_8 
	\end{pmatrix} \begin{pmatrix}
	\partial_p \mathcal{H} \\
	\partial_q \mathcal{H}
\end{pmatrix},\label{ec:FullPHS_PDE}
\end{equation}
where
\[
\begin{aligned}
	\mathcal{J} &= \mathcal{P}_1\frac{\partial(\cdot)}{\partial \xi_1}+\mathcal{P}_2\frac{\partial(\cdot)}{\partial \xi_2}+\mathcal{P}_0,\\
	-\mathcal{J}^* &= \mathcal{P}_1^T\frac{\partial(\cdot)}{\partial \xi_1}+\mathcal{P}_2^T\frac{\partial(\cdot)}{\partial \xi_2}-\mathcal{P}_0^T,\\
	\partial_p \mathcal{H} &= \frac{\partial \mathcal{H}}{\partial p}, \hspace{40pt}
	\partial_q \mathcal{H} = \frac{\partial \mathcal{H}}{\partial q}, \\
	\mathcal{H}(t,\xi_1,\xi_2) &= \frac{1}{2} q^TC_q q+ p^TM_p^{-1}p.
\end{aligned}
\]
The matrices $\mathcal{P}_1$, $\mathcal{P}_2$ and $\mathcal{P}_0$ are given in Appendix A.

This latter function is referred to as the energy density, since the Hamiltonian can be written as 
\[
	H(t) = \iint_\Omega \mathcal{H}(t,\xi_1,\xi_2)dA.
\]
The time derivative of the total energy is then given by
\[
\dot{H} = \oint_{\partial \Omega} \begin{pmatrix}
	(\partial_p \mathcal{H})^T\mathcal{P}_1 (\partial_q \mathcal{H}) & (\partial_p \mathcal{H})^T\mathcal{P}_2 (\partial_q \mathcal{H})
\end{pmatrix}\cdot \hat{n}\; ds,
\]
where $\hat{n}$ is the outward-pointing unit normal vector on the boundary $\partial \Omega$. With this in mind, the boundary inputs and outputs can be defined as linear combinations of the generalized velocities $\partial_p \mathcal{H}$ and the generalized stresses $\partial_q \mathcal{H}$.

\subsection{Problem of interest and discretization}\label{subsec:Discretization}
\begin{figure}[ht]
	\centering
	\includegraphics[width=.5\linewidth]{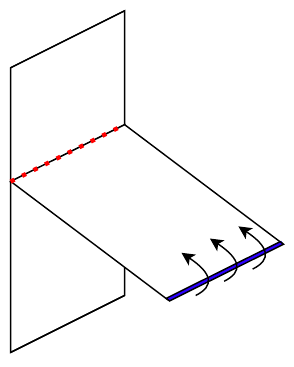}
	\caption{Example of problem of interest.}
	\label{fig:exampPoI}
\end{figure}

\begin{figure}[ht]
	\centering
	\includegraphics[width=.8\linewidth]{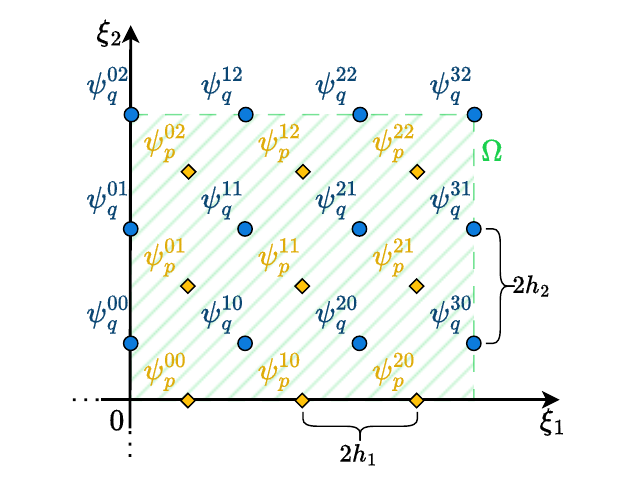}
	\caption{Example of problem of interest.}
	\label{fig:exampPoI_disc}
\end{figure}

For the problem of interest, we consider a plate that is clamped at $\xi_2=0$, the boundaries at $\xi_1=0$ and $\xi_2=L_1$ are free; therefore, no external forces are applied on these edges. Finally, the boundary at $\xi_2 = L_2$ is subject to force inputs, where the finite-dimensional controller will be designed and interconnected. The configuration of the problem is illustrated in Fig.~\ref{fig:exampPoI}.
To design the controller, we first discretize the plate model using a centered finite-difference scheme on staggered grids to preserve its port-Hamiltonian structure, as presented in \citep{Diaz_ECC2026}. To do so, we define the sets of points $\Psi_q$, $\Psi_p$ which are subsets of $\Omega$ given by
\begin{align*}
	\Psi_q&:\Big\{\psi_q^{mn}=\Big(h_1\cdot(2m+1)\;,\;h_2\cdot(2n)\Big)\Big|\psi_q^{mn}\in\Omega\Big\},\\
	\Psi_p&:\Big\{\psi_p^{mn}=\Big(h_1\cdot(2m)\;,\;h_2\cdot(2n+1)\Big)\Big|\psi_p^{mn}\in\Omega\Big\},
\end{align*}
with the pair $(m,n)$ two positive integers. This can be exemplified by Figure \ref{fig:exampPoI_disc}. 
In this example, the points $\psi_p^{\{0,1,2\}0}$ at the boundaries correspond to the clamped part; the points $\psi_q^{0\{0,1\}}$ and $\psi_q^{3\{0,1\}}$ at the boundaries correspond to the free borders; and finally the points $\psi_q^{\{0,1,2,3\}2}$ correspond to the boundaries where our controller will act.

Then, we consider the set of discrete general coordinates
\begin{align*}
	\mathcal{X}_q&:\Big\{x_q^{mn}=q(t,\psi_q^{mn})|\;\forall \psi_q^{mn}\in\Psi_q\setminus \partial \Omega\Big\} ,\\
	\mathcal{X}_p&:\Big\{x_p^{mn}=p(t,\psi_p^{mn})|\;\forall \psi_p^{mn}\in\Psi_p\setminus \partial \Omega\Big\} .
\end{align*}
Consider an ordered vector $x_q$ in $\mathcal{X}_q$ where \linebreak ${\forall(m,n)\psi_q^{mn}\in\Psi_q\setminus \partial \Omega}$ and similarly for $x_p$ in $\mathcal{X}_p$, we can write the discretized Hamiltonian as
\begin{equation}
	H_d(x_q,x_p)=\frac{1}{2}\Big(x_q^TC_dx_q+x_p^TM_d^{-1}x_p\Big)4h_1h_2,
	\label{ec:2DHamiltonian_Discrete}
\end{equation}
where $M_d$ is a block diagonal matrix including $M_p(\xi_1,\xi_2)$ evaluated at each $\psi_p^{mn}\in\Psi_p\setminus \partial \Omega$ in the same order of $x_p$. Similarly $C_d$ is a block diagonal matrix including $C_q(\xi_1,\xi_2)$ evaluated at each $\psi_q^{mn}\in\Psi_q\setminus \partial \Omega$ in the same order of $x_p$. We can also write the dynamic equations 
\begin{equation}
	\begin{aligned}
		\begin{pmatrix}
			x_{p}(t)\\
			x_{q}(t)
		\end{pmatrix} &= \begin{pmatrix}
			\mathbf{0} & \mathcal{P}_{dp} \\
			\mathcal{P}_{dq} & \mathbf{0}
		\end{pmatrix}\begin{pmatrix}
			M_d^{-1}x_p(t) \\
			C_dx_q(t)
		\end{pmatrix}+B_d \begin{pmatrix}
			u_q(t)\\
			u_p(t)
		\end{pmatrix} \\
		\begin{pmatrix}
			y_p(t)\\
			y_q(T)
		\end{pmatrix} &= {B_d}^T\begin{pmatrix}
			M_d^{-1}x_p(t) \\
			C_dx_q(t)
		\end{pmatrix}
	\end{aligned}\label{ec:discretised_ODE_2D}
\end{equation}
where 
\[
\begin{aligned}
	\mathcal{P}_{dp}&=\mathcal{I}_1^{(p)}\otimes\mathcal{P}_1h_2+\mathcal{I}_2^{(p)}\otimes\mathcal{P}_2h_1+\mathcal{I}_0^{(p)}\otimes\mathcal{P}_0h_1h_2,\\
	\mathcal{P}_{dq}&=\mathcal{I}_1^{(q)}\otimes{\mathcal{P}_1}^Th_2+\mathcal{I}_2^{(q)}\otimes{\mathcal{P}_2}^Th_1-\mathcal{I}_0^{(q)}\otimes{\mathcal{P}_0}^Th_1h_2,
\end{aligned}
\]
with $\mathcal{I}_i^{(j)}$ representing the connection coefficients related to the $\mathcal{P}_i$ matrix for the points in the $\Psi_j$ set. And where
\[
B_d = \frac{1}{4h_1h_2}\begin{pmatrix}
	B_p & \mathbf{0} \\
	\mathbf{0} & B_q 
\end{pmatrix},
\]
with
\[
B_p = \mathcal{I}^{(up)}\otimes I_n,
\]
\[
B_q = \left(\mathcal{I}_1^{(u)}\otimes{\mathcal{P}_1}^Th_2+\mathcal{I}_2^{(u)}\otimes{\mathcal{P}_2}^Th_1-\mathcal{I}_0^{(u)}\otimes{\mathcal{P}_0}^Th_1h_2\right),
\]
where $\mathcal{I}^{(up)}$ represents the mapping from $\Psi_q \cap \partial \Omega$ to $\Psi_p$, $n$ is the dimension of the generalized coordinate $p$, and $\mathcal{I}_i^{(u)}$ represents the mapping from $\Psi_p \cap \partial \Omega$ to $\Psi_q$. It can be easily checked that $P_{dp}=-{P_{dq}}^T$, and from this \eqref{ec:discretised_ODE_2D} is the approximation of \eqref{ec:FullPHS_PDE} in the PH structure that we will use for the rest of this paper.

\section{Controller Design}
This section recalls the control design strategy proposed in \citep{toledo2020}, then addresses the controllability issues introduced by the spatial discretization. Finally, the implementation procedure of the controller for the Mindlin plate model is presented.
\subsection{Observer based state feedback for 1D PHS}\label{sec:Controldesign}
In \citep{toledo2020}, an early lumped observer-based state-feedback (OBSF) control design method is proposed for one-dimensional boundary controlled port-Hamiltonian systems. The method consists of designing the state-feedback gain $K$ and the observer gain $L$ for a finite-dimensional approximation of the distributed port-Hamiltonian system such that the resulting OBSF controller is strictly positive real, i.e., it admits a port-Hamiltonian representation with a strictly positive controller dissipation matrix $R_c$. Under this condition, asymptotic closed-loop stability is guaranteed when the finite-dimensional OBSF controller is interconnected with the distributed port-Hamiltonian system, thereby overcoming the spillover effect.

The procedure can be summarized by the following five-step algorithm.
\begin{enumerate}
	\item[(\refstepcounter{enumi}\arabic{enumi}.A)\label{alg:Disc}] Derive a finite-dimensional approximation of the boundary control problem in PH structure, given by
	\begin{equation}
		\begin{aligned}
			\dot{x}_d(t) &= (J_d-R_d)Q_dx_d(t)+B_du(t),\\
			y(t) &= {B_d}^TQ_dx_d(t).
		\end{aligned}\label{ec:discretePHS}
	\end{equation}
	\item[(\refstepcounter{enumi}\arabic{enumi}.A)\label{alg:Kdesign}] Design $K$ such as $A_K = (J_d-R_d)Q_d-B_dK$ for all eigenvalues $\operatorname{Re}\{\lambda_i\}<0$. 
	\item[(\refstepcounter{enumi}\arabic{enumi}.A)\label{alg:Hm}] Choose a positive definite matrix $R_c$ such that
	\begin{equation}
		H_m = \begin{pmatrix}
			A_K & 2R_c\\ -C_K & -{A_K}^T
		\end{pmatrix}\label{ec:Hm_forARE}
	\end{equation}
    with
    \[
    C_K = -(K^T{B_d}^TQ_d+Q_d^TB_dK)
    \]
	has no pure imaginary eigenvalues.
	\item[(\refstepcounter{enumi}\arabic{enumi}.A)\label{alg:ARE}] Solve the algebraic Ricatti equation (ARE) given by
	\begin{equation}
		{A_K}^TQ_c+Q_cA_K+2Q_cR_cQ_c+C_K=0,
		\label{ec:ARE}
	\end{equation}
	for the matrix $Q_c$.
	\item[(\refstepcounter{enumi}\arabic{enumi}.A)\label{alg:OBSF}] Implement the OBSF controller 
	\begin{equation}\label{eq:OBSFcontroller}
		\begin{aligned}
			\dot{x}_c(t) &=  (J_c-R_c)Q_cx_c(t)+B_cu_c(t)+B_d r(t) \\
            y_c(t) &= B_c^TQ_c x_c(t) \\
			u(t) &= -y_c(t)+r(t) \\
            u_c(t) &= y(t)
		\end{aligned}
	\end{equation}
	with $B_c=L={Q_c}^{-1}K^T$.
\end{enumerate}

\subsection{Controllable mode projection}
In the design procedure (\ref{alg:Disc}.A), the discretized system \eqref{ec:discretePHS} is assumed to be controllable and observable, whereas in \eqref{ec:discretised_ODE_2D} the matrices $P_{dp}$ and $P_{dq}$ are rectangular. As a consequence, the interconnection matrix exhibits a zero block associated with linear combinations of strain approximations, which can be interpreted as geometric constraints introduced by the discretization. These constraints may, in turn, give rise to uncontrollable modes.
To identify these uncontrollable modes in the discretized system, we exploit the fact that every skew-Hermitian (and therefore every skew-symmetric) matrix is similar to a block diagonal matrix (\citep{becker1973}), i.e.
\[
	U^TJU=D
\]
with $U$ a unitary matrix and
\[
	D = \begin{pmatrix}
		\text{diag}\Big(\begin{bmatrix}
			0 & -\sqrt{\alpha_1} \\ \sqrt{\alpha_1} & 0
		\end{bmatrix},\;\dots,\;\begin{bmatrix}
			0 & -\sqrt{\alpha_k} \\ \sqrt{\alpha_k} & 0
		\end{bmatrix}\Big) & \mathbf{0}\\ \mathbf{0} & \mathbf{0}
	\end{pmatrix},
\]
where $\mathbf{0}$ denotes a zero matrix of appropriate dimension, \linebreak${\alpha_i\,|\,\forall i\in \{1,\dots,k\}}$ are the nonzero eigenvalues of $J^TJ$.
\begin{remark}
 The eigenvalues $\alpha_i$ are defined by the norm of the eigenvalues of $J$ and therefore always have even algebraic multiplicity, defined by $2\ell$. Consequently, for each $\alpha_i$ appearing in $D$, the corresponding $2 \times 2$ block occurs exactly $\ell$ times. In particular, if $\alpha_i$ has algebraic multiplicity $2$, $D$ only has one block dependant of $\alpha_i$, whereas if it has algebraic multiplicity $4$, it appears in two blocks of $D$.

\end{remark}
A change of coordinates $x_u=U^Tx$ can be applied to the port-Hamiltonian system
\begin{center}
	\[\begin{aligned}
		\dot{x}(t) &= JHx(t)+Bu(t)\\
		y(t) &= B^THx(t)
	\end{aligned}\to\begin{aligned}
		\dot{x}_u(t) &= D \bar{H}x_u(t)+\bar{B}u(t)\\
		y(t) &= \bar{B}^T\bar{H}x_u(t)
	\end{aligned},\]
\end{center}
where $\bar{H}=U^THU$, $\bar{B}=U^TB$. Moreover, the Hamiltonian in the new coordinates is given by
\[
x^THx = x^TUU^THUU^Tx = x_u^T\bar{H}x_u.
\]
In the new coordinates, the system can always be decomposed as
\[\begin{aligned}
	\dot{x}_c(t) &= D_c H_c x_c(t)+B_cu(t),\\
	\dot{x}_{\bar{c}}(t) &= B_{\bar{c}}u(t),
\end{aligned}\]
where $x_c$ contains the first $2k$ variables of $x_u$, $x_{\bar{c}}$ the remaining $n-2k$ variables. Furthermore, $H_c$ denotes the upper-left $2k$ block of $\bar{H}$, $B_c$ the upper $2k\times n_u$ block of $\bar{B}$, and $B_{\bar{c}}$ the lower $(n-2k)\times n_u$ matrix of $\bar{B}$ and
\[
D_c = \text{diag}\Big(\begin{bmatrix}
	0 & -\sqrt{\alpha_1} \\ \sqrt{\alpha_1} & 0
\end{bmatrix},\;\dots,\;\begin{bmatrix}
	0 & -\sqrt{\alpha_k} \\ \sqrt{\alpha_k} & 0
\end{bmatrix}\Big).
\]
Therefore, the states $x_{\bar{c}}(t)$ have no internal dynamics and do not affect the evolution of $x_c(t)$. In practice, the matrix $B_{\bar{c}}$ is usually a zero matrix, and these states represent geometric constraints induced by the discretization. For simplicity, assume that the plate is undeformed at an initial time $t_0$ and that $B_{\bar{c}} \equiv 0$. Then,
\[
	\int_{t_0}^{t}\dot{x}_{\bar{c}}(\tau)d\tau=0
\]
which implies $x_{\bar{c}}(t)=0$. Hence, the finite-dimensional port-Hamiltonian system can be expressed in terms of the controllable states as
\begin{equation}
	\begin{aligned}
		\dot{x}_c(t)&=D_c\bar{H}_cx_c(t)+B_cu(t)\\
		y(t) &= B_c^T\bar{H}_cx_c(t)
	\end{aligned},
\end{equation}
and the OBSF design algorithm presented in Subsection \ref{sec:Controldesign} can therefore be applied.

\subsection{OBSF controller design for Mindlin plate\label{ssec:OBSF}}
To employ the previously mentioned OBSF design algorithm, the state feedback gain $K$ is first designed as indicated in (\ref{alg:Kdesign}.A). 
To this end, the Linear Quadratic Regulator (LQR) method is adopted with the cost functional
\[
J_{LQR}(u)=\int_0^\infty \begin{bmatrix}
	x^T & u^T
\end{bmatrix}\underbrace{\begin{bmatrix}
		Q & N\\ N^T & R
\end{bmatrix}}_{\mathcal{J}_{LQR}}\begin{bmatrix}
	x\\u
\end{bmatrix}dt,
\]
where the weighting matrix $\mathcal{J}_{LQR}$ is required to be semi positive definite. Since the energy of the system is proportional to $x_c^{T}\bar{H}_c x_c$ and the power supplied through the inputs depends on $x_c^{T}\bar{H}_c B_c u$, both quantities can be simultaneously minimized by selecting 
\[Q=\beta^2\bar{H}_c, \;\; N=\bar{H}_cB_c/2, \;\; R=B_c^T\bar{H}_cB_c/(4\beta^2)\] 
where $\beta>0$is a tuning parameter that balances the relative weight between the stored energy and the input power. For numerical stability, additional positive semi-definite terms may be added to $Q$ and $R$.

\begin{figure}[ht]
	\centering
	\includegraphics[width=0.85\linewidth]{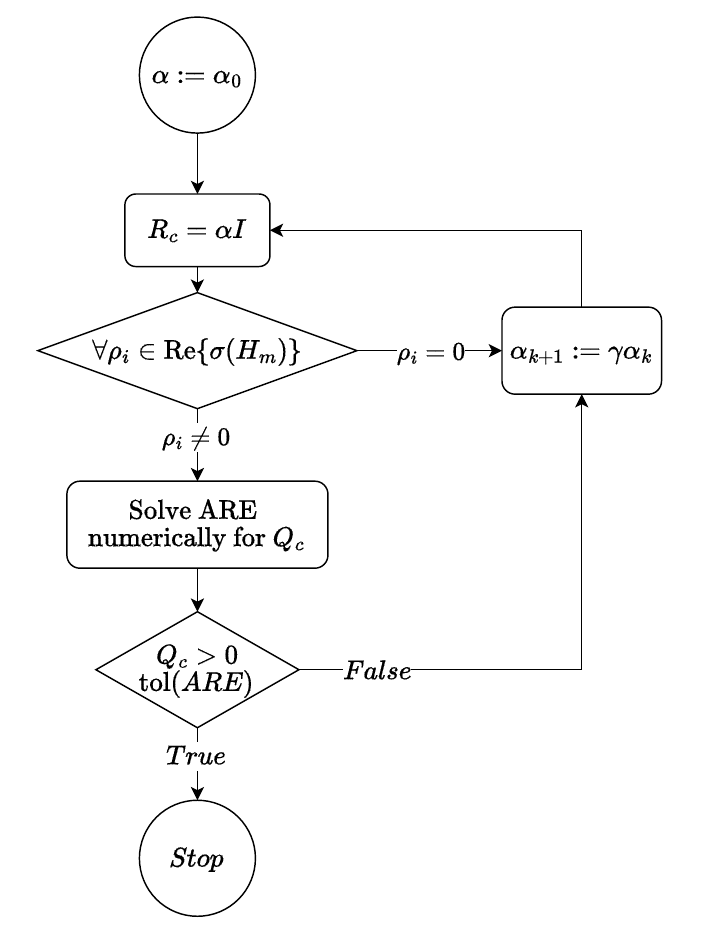}
	\caption{$R_c$ matrix calculation.}
	\label{fig:AlgorithmRc}
\end{figure}

Once the feedback gain $K$ has been obtained via the LQR method, a positive definite matrix $R_c$ must be determined such that the matrix \eqref{ec:Hm_forARE} has no purely imaginary eigenvalues. This guarantees that the algebraic Riccati equation \eqref{ec:ARE} in (\ref{alg:ARE}.A) admits a symmetric positive definite solution $Q_c = Q_c^{T} > 0$. To this end, $R_c$ is chosen as $R_c=\alpha I_{n_c}$ where $n_c = \dim(x_c)$. The scalar $\alpha$ is determined iteratively as follows. Starting from an initial value of $\alpha$, the spectrum $\sigma(H_M)$ is evaluated. If an eigenvalue lies on the imaginary axis (within a prescribed tolerance), $\alpha$ is reduced by a factor $\gamma \in (0,1)$, i.e., $\alpha \leftarrow \gamma \alpha$. The Riccati equation \eqref{ec:ARE} is then solved, and the positivity of $Q_c$ as well as the residual of the ARE are checked. If either condition is not satisfied, $\alpha$ is further reduced by the same factor. This procedure is summarized in Fig.~\ref{fig:AlgorithmRc}.

Once the Riccati equation \eqref{ec:ARE} is solved and a symmetric positive definite solution $Q_c = Q_c^{T} > 0$ is obtained, the OBSF controller is computed from \eqref{eq:OBSFcontroller}.

\section{Numerical Results}
This section presents a numerical implementation of the algorithm developed in the previous section and follows them with some comments about the spillover effect.
The simulation considers an aluminum plate of thickness $10 \mathrm{mm}$, width $1 \mathrm{m}$, and length $2 \mathrm{m}$, with physical parameters $E = 68\mathrm{GPa}$, $\nu = 0.36$, $\rho = 2698.9\mathrm{kg/m^3}$, $G = 25\mathrm{GPa}$, and $\kappa = \pi^2/12$.

\subsection{Controller Implementation}
\begin{figure}[ht]
    \centering
    \includegraphics[width=.95\linewidth]{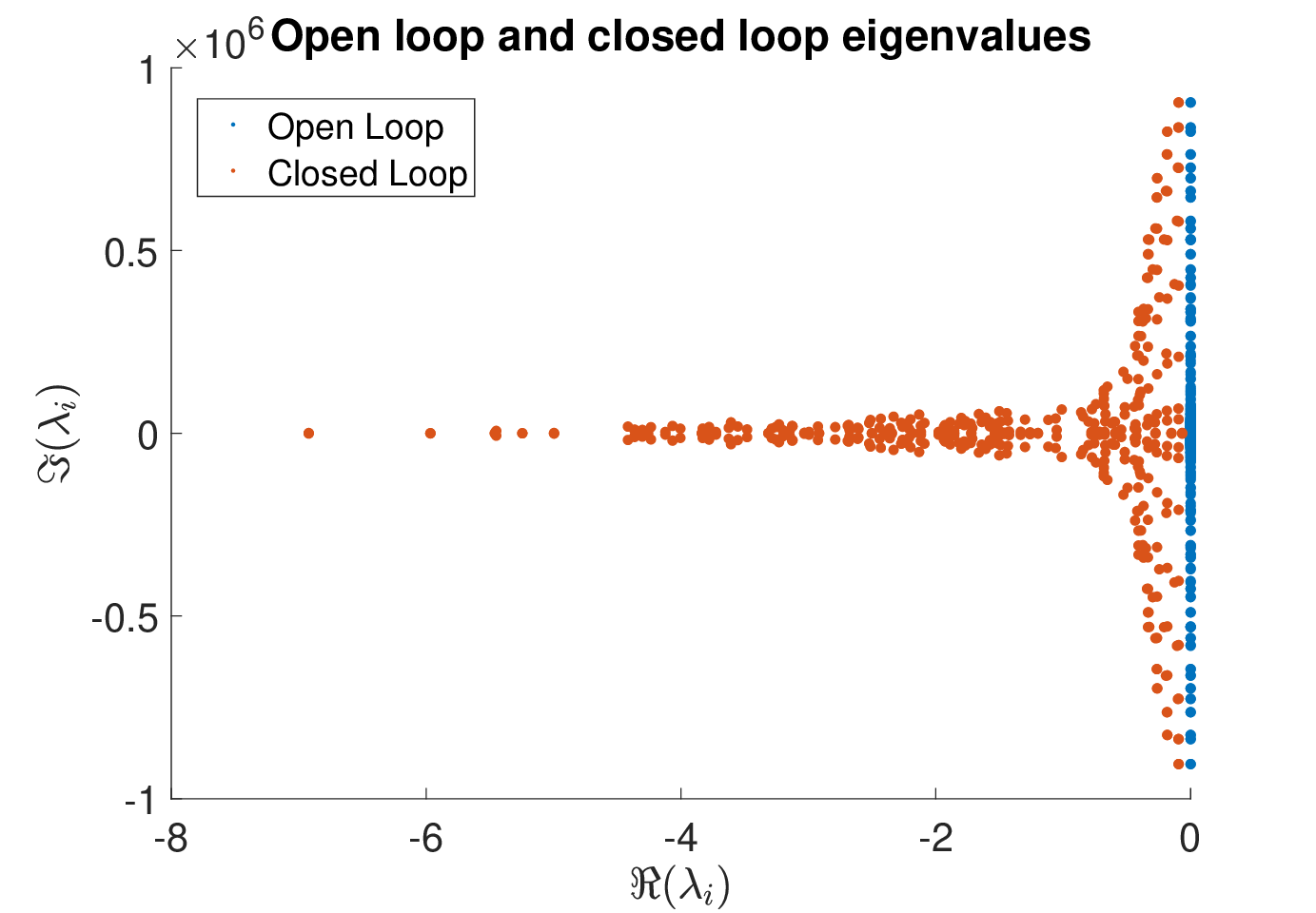}
    \caption{Eigenvalues of different matrices.}
    \label{fig:eigenvalues6x6}
\end{figure}

First, a low-order spatial discretization of the plate is obtained using the method presented in Subsection \ref{subsec:Discretization}, with $42$ internal points for the kinetic energy domain and $36$ internal points for the potential energy domain. This results in a discrete state-space model of dimension $498$, which can be reduced to a controllable subspace of dimension $420$. Then, to design the controller, we apply the LQR-based control design described in Section \ref{ssec:OBSF} to compute the state-feedback gain $K$. Using this gain, a stable closed-loop system is obtained under full-state feedback. The corresponding eigenvalues are shown in Fig. \ref{fig:eigenvalues6x6}. The red markers denote the eigenvalues of the open-loop plate model, which all lie on the imaginary axis due to the absence of damping. The blue markers represent the eigenvalues of the closed-loop system with state feedback, all of which are located in the left half of the complex plane.

Furthermore, we take $R_c=\alpha I_{n_c}$ with $\alpha=0.0114$ which results in the eigenvalues of the matrix in \eqref{ec:Hm_forARE} having real parts close to $1$. This choice enables us to solve the algebraic Riccati equation (ARE) in (\ref{alg:ARE}.A) using the \verb|icare| command from the \verb|Control System Toolbox| in MATLAB. Since the $C_K$ matrix is sign indefinite, the ARE must be reformulated in its inverse form,
\[
    Q_c^{-1}{A_K}^T+A_KQ_c^{-1}+2R_c+Q_c^{-1}C_KQ_c^{-1}=0,
\]
where $2R_c$ is positive definite by construction. Accordingly, the command \verb|icare(Ak',[],2*Rc,[],[],[],Ck)| is implemented. Due to numerical sensitivity, it is essential to verify that the resulting matrix $Q_c$ is effectively positive definite. Finally, the passive OBSF can be derived through the formula \eqref{eq:OBSFcontroller} with  $Q_c$. 
To evaluate the performance of the proposed controllers, we consider a consistent reference signal
\[
r(t) = u_{ref}\,f_{act}(t),
\]
where $f_{act}(t)\in[-1,1]$ is a modulating function over time.

\subsection{Simulation results }
For simulation purposes, a midpoint time-discretization scheme is employed, as it preserves the Hamiltonian of the system. To assess the performance and robustness of the proposed controller, three observer-based state-feedback (OBSF) strategies are considered for stabilizing the Mindlin plate model: a perfect OBSF, an arbitrary Luenberger OBSF, and the passive OBSF designed in this work. The perfect OBSF assumes full access to the system states at all times, i.e., full-state feedback. The same reference signal is applied in all three cases, with the controller activated at $t=2\,s$. The reference is turned off by setting $f_{act}$ at $t=6\,s$ and subsequently switched to $f_{act}=-0.75$ for $t>=10\,s$. 
\begin{figure}[ht]
    \centering
    \includegraphics[width=.95\linewidth]{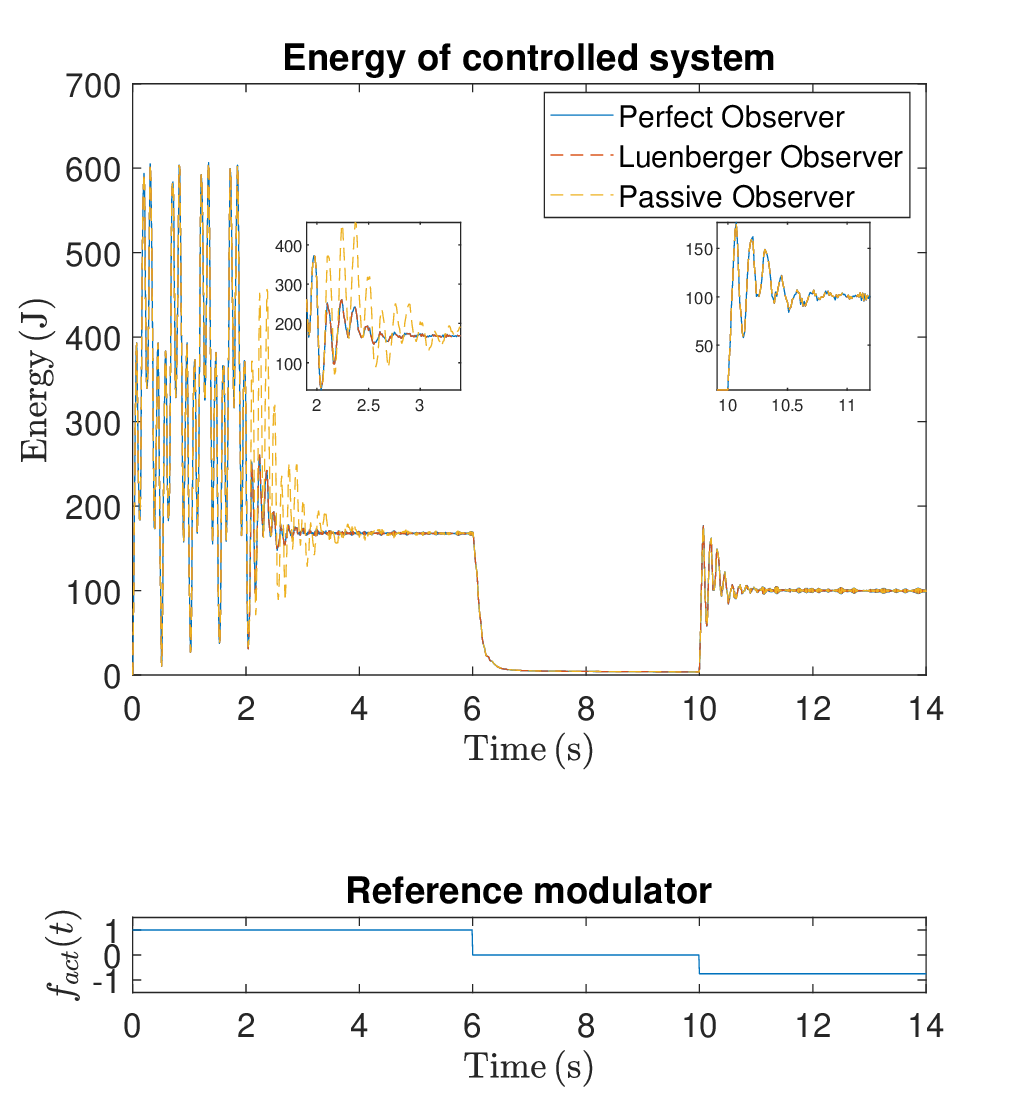}
    \caption{Energy overtime for controlled systems}
    \label{fig:HamiltonianComparison}
\end{figure}
To compare the performance of the different systems, the Hamiltonian evolution is shown in Fig. \ref{fig:HamiltonianComparison}. The blue solid line corresponds to the closed-loop system with the perfect OBSF (full-state feedback), while the yellow dashed line represents the Luenberger OBSF case; both exhibit nearly identical stabilization times. The orange solid line shows the Hamiltonian of the closed-loop system with the proposed passive OBSF, which converges to the same equilibrium with a longer convergence time of approximately 2 seconds. This behavior arises because the passive OBSF preserves the Hamiltonian structure, whereas the observer gain cannot be chosen freely; this restriction represents a necessary compromise to guarantee the stability of the infinite-dimensional system and to avoid spillover effects. We next demonstrate that the passive OBSF is able to stabilize the infinite-dimensional plate system, whereas the Luenberger OBSF does not ensure stability in this case.

\begin{figure}[ht]
    \centering
    \includegraphics[width=.95\linewidth]{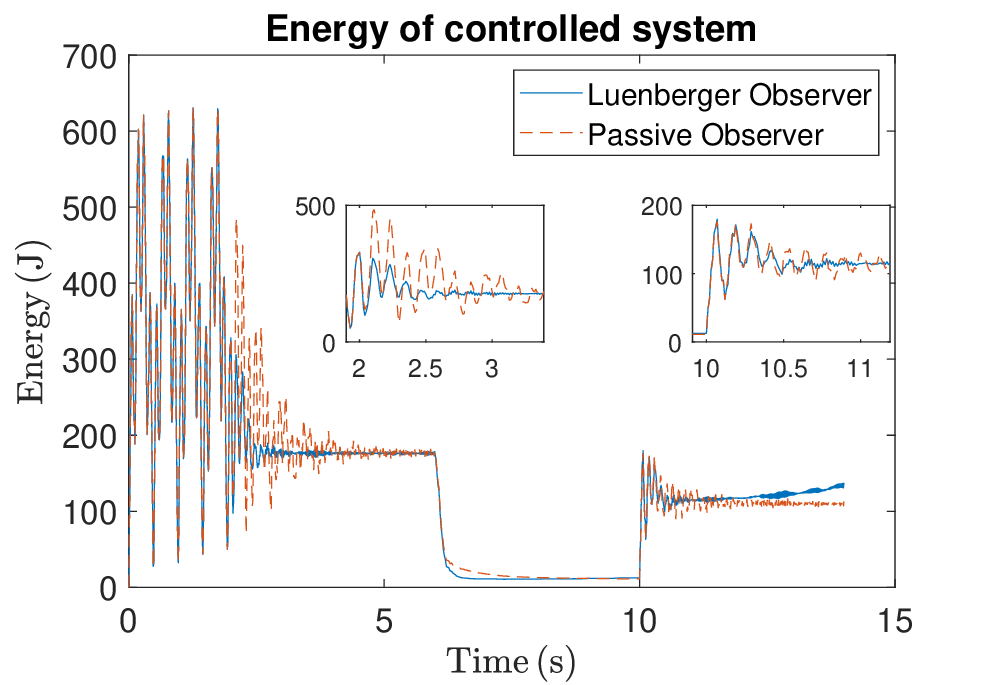}
    \caption{Energy overtime for controlled systems}
    \label{fig:HamiltonianComparison_HO}
\end{figure}

To this end, a higher-order spatial discretization is employed, using 420 internal points for the kinetic energy and 400 internal points for the potential energy to approximate the infinite-dimensional plate. For the interconnection at the boundary, we exploit the fact that the velocities follow a linear interpolation between neighbouring points. This leads to express $v_1 = \Phi v_2$, where $v_1$ are the velocities of the low order system, $v_2$ are the velocities of the high order system and $\Phi$ is a matrix that comes from the linear combination given by the interpolation. The forces must be defined such as the interconnection satisfies a power-preserving connection and therefore $F_2 = \Phi^T F_1$. The Hamiltonian evolutions of both controlled models are shown in Fig. \ref{fig:HamiltonianComparison_HO}.

\begin{figure}[ht]
    \centering
    \includegraphics[width=.8\linewidth]{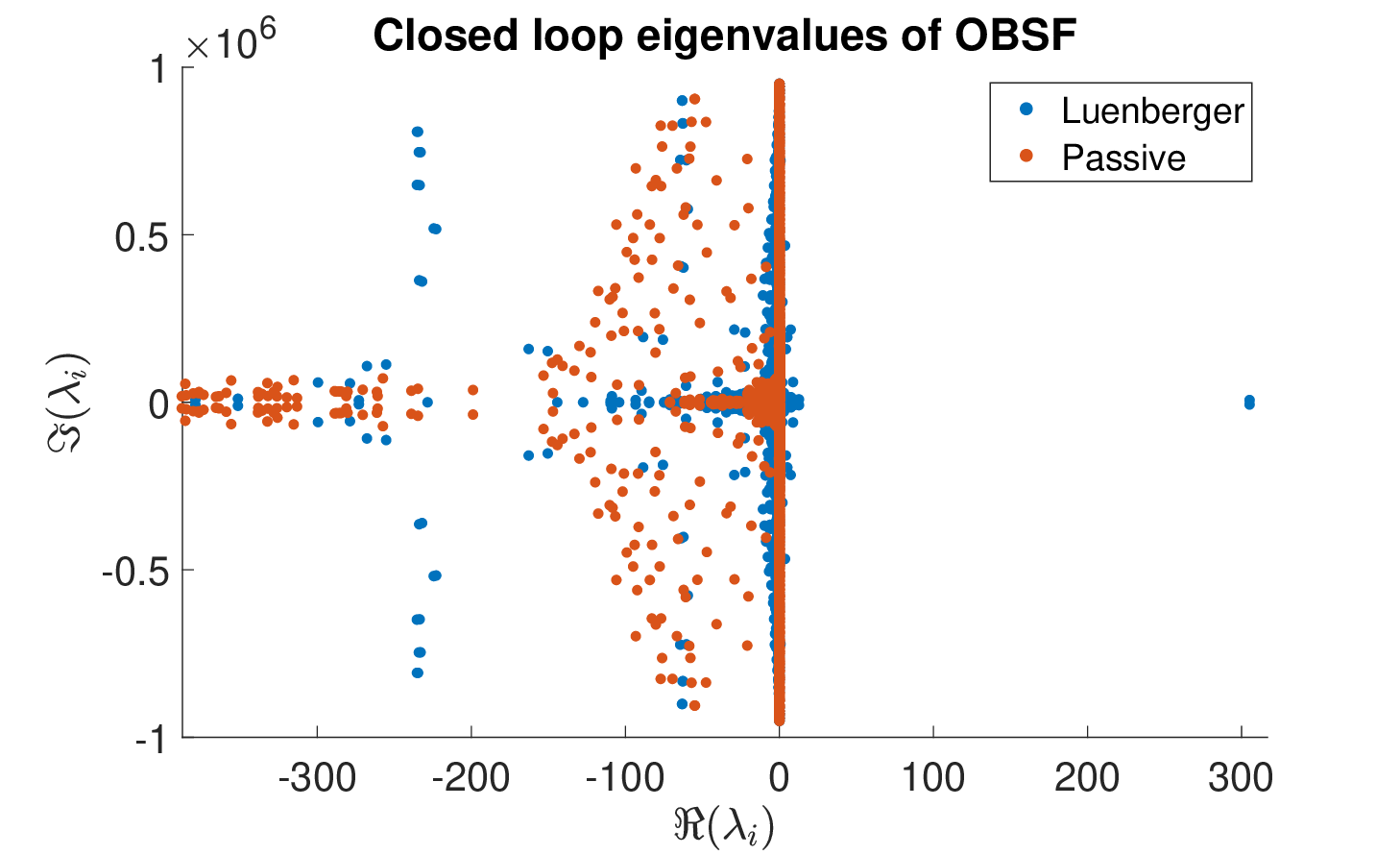}
    \caption{Eigenvalues for closed loop systems}
    \label{fig:eigenvalues_HO}
\end{figure}

It can be observed that the energy of the system controlled by the Luenberger OBSF begins to diverge toward the end of the simulation. This behavior is attributed to the spillover effect, as also reported in Fig. 1 of \cite{toledo2020}. This observation is further confirmed in Fig. \ref{fig:eigenvalues_HO}, where several eigenvalues of the closed-loop system with the Luenberger OBSF (blue dots) are located in the right half-plane, leading to instability. In contrast, all eigenvalues of the closed-loop system with the passive OBSF (blue markers) remain in the left half-plane, thereby guaranteeing stability.

\section{Conclusion}
In this paper, an observer-based state-feedback (OBSF) control strategy for a boundary-controlled two-dimensional port-Hamiltonian Mindlin plate is investigated using an early lumping approach. First, a structure-preserving discretization based on staggered grids is employed to obtain a finite-dimensional approximation of the 2D Mindlin plate. The resulting discretized model is then decomposed into a controllable subsystem and a set of geometric constraints induced by the discretization. Subsequently, the state-feedback and observer gains are designed using a low-order finite-dimensional approximation, ensuring that the resulting OBSF controller admits a fully damped port-Hamiltonian realization (strictly positive real). Finally, numerical simulations demonstrate that the closed-loop system remains stable even when the finite-dimensional OBSF controller is interconnected with a higher-order 2D Mindlin plate model used in simulation.
\\Future research will address the well-posedness of the closed-loop system and the experimental validation of the proposed control strategy.

\bibliography{bibliography}

\appendix
\section{System matrices}
The matrices of the interconnection operator is given by
\begin{align}
    \mathcal{P}_1 &= \begin{pmatrix}
        1 & 0 & 0 & 0 & 0 & 0 & 0 & 0\\
        0 & 0 & 1 & 0 & 0 & 0 & 0 & 0\\
        0 & 0 & 0 & 0 & 0 & 0 & 1 & 0\\
        0 & 0 & 0 & 1 & 0 & 0 & 0 & 0\\
        0 & 0 & 0 & 0 & 0 & 1 & 0 & 0
    \end{pmatrix},\\
    \mathcal{P}_2 &= \begin{pmatrix}
        0 & 0 & 1 & 0 & 0 & 0 & 0 & 0\\
        0 & 1 & 0 & 0 & 0 & 0 & 0 & 0\\
        0 & 0 & 0 & 0 & 0 & 0 & 0 & 1\\
        0 & 0 & 0 & 0 & 0 & 1 & 0 & 0\\
        0 & 0 & 0 & 0 & 1 & 0 & 0 & 0
    \end{pmatrix},\\
    \mathcal{P}_0 &= \begin{pmatrix}
        0 & 0 & 0 & 0 & 0 & 0 & 0 & 0\\
        0 & 0 & 0 & 0 & 0 & 0 & 0 & 0\\
        0 & 0 & 0 & 0 & 0 & 0 & 0 & 0\\
        0 & 0 & 0 & 0 & 0 & 0 & 1 & 0\\
        0 & 0 & 0 & 0 & 0 & 0 & 0 & 1
    \end{pmatrix}.
\end{align}
\end{document}